\documentclass[11pt,letterpaper]{article}
\usepackage[margin=1in]{geometry}
\usepackage[T1]{fontenc}
\usepackage{lmodern}
\usepackage{amsmath,amssymb,amsthm,mathtools}
\usepackage{microtype}
\usepackage[hidelinks]{hyperref}
\allowdisplaybreaks
\newcommand{\doi}[1]{DOI: \href{https://doi.org/#1}{\nolinkurl{#1}}}
\newtheorem{theorem}{Theorem}
\theoremstyle{remark}
\newtheorem{remark}{Remark}
\theoremstyle{plain}

\newtheorem{lemma}{Lemma}
\newtheorem{corollary}{Corollary}
\newtheorem{proposition}{Proposition}
\newtheorem{claim}{Claim}
\newenvironment{claimproof}{%
  \begin{proof}[Proof of the claim]}%
  {\end{proof}}

\title{The Common-Color Spectrum of\\
Endpoint Set-Coloring Ramsey Numbers}
\author{
Guantao Chen\thanks{\protect\raggedright Department of Mathematics and Statistics, Georgia State University,
Atlanta, GA 30303, USA. E-mail: \texttt{gchen@gsu.edu}.}
\and
Qiong Fan\thanks{\protect\raggedright School of Mathematics and Statistics, Central China Normal University,
Wuhan, China. E-mail: \texttt{xfanqiong@ccnu.edu.cn}.}
\and
Chunxiang Wang\thanks{\protect\raggedright Corresponding author. School of Mathematics and Statistics,
Central China Normal University, Wuhan, China.
E-mail: \texttt{wcxiang@mail.ccnu.edu.cn}.}
}
\date{}

\begin{document}
\maketitle

\begin{abstract}
The set-coloring Ramsey number $R(n;r,s,t)$ is the least $N$
such that every assignment of an $s$-element subset of
$[r]=\{1,\ldots,r\}$ to each edge of $K_N$ contains a $K_n$
whose edge labels have at least $t$ common colors.
Motivated by the work of Erd\H{o}s and Szemer\'edi on
$R(n;r,r-1,1)$, we write
$R(n;r,r-1,t)=\exp\bigl(nf_{r,n}(t)\bigr)$ and study
$f(t):=f_{r,n}(t)$ as $t$ ranges from $1$ to $r-1$.

For fixed integers $r\ge2$ and $1\le t<r$, we prove that,
as $n\to\infty$,
\(
 \frac12\log\frac{r}{r-t}+o(1)
 \le f(t)
 \le \frac{r}{r-t}\log\frac{r}{r-t}
     +C+O_{r,t}\!\left(\frac{\log n}{n}\right),
\)
where $C>0$ is an absolute constant.
Comparisons between different numbers of colors show that
the rates $\liminf_{n\to\infty}f_{r,n}(t)$ and
$\limsup_{n\to\infty}f_{r,n}(t)$, defined for fixed $r$ and $t$,
each have a limit as $r\to\infty$ and
$\frac{t}{r}\to\alpha\in(0,1)$ from above.
The resulting two functions of $\alpha$ are finite, positive,
nondecreasing, and right-continuous.
Writing $\Delta_t=f(t+1)-f(t)$, we further prove that,
as $r\to\infty$, uniformly for $n\ge r^2$, at least
$\frac{(1-o(1))\log r}{\log\log r}$ indices
$1\le t<\frac{r}{2}$
satisfy $\Delta_t\ge\frac{1-o(1)}{2r}$.
For $r,n\ge3$ and $2\le h\le r-1$, the final $h-1$
increments account for at least
$1-\frac{1}{\left\lceil h/2\right\rceil}$
of the total increase $f(r-1)-f(1)$.
For each fixed $r\ge4$, the last increment has a positive
lower limit as $n\to\infty$.
Moreover, for $r,n\ge3$, it exceeds the sum of all preceding
increments whenever $r$ is even or $n=3$.

\end{abstract}


\section{Introduction}
All graphs are finite, simple, and undirected.
Write $[r]=\{1,\ldots,r\}$, and let $K_N$ denote the complete
graph on $N$ vertices. For a graph $G$, write $E(G)$ for its
edge set, $e(G)=|E(G)|$ for its number of edges, and
$\overline G$ for its complement.
All logarithms are natural.
Unless stated otherwise, bounds in $n$ are understood
for fixed color parameters and all sufficiently large $n$;
the threshold may depend on those parameters.

For integers $n\ge3$ and $r\ge2$, the multicolor Ramsey number
$R(n;r)$ is the least $N$ such that every edge-coloring
$\chi:E(K_N)\to[r]$ contains a monochromatic $K_n$.
An $r$-coloring may use fewer than $r$ colors, and $R(n;1)=n$.
The classical bounds are
\[
 2^{cnr}\le R(n;r)\le2^{c'nr\log r},
\]
where $c,c'>0$ are absolute constants;
see~\cite{CFHMSV2024}.

Set-coloring Ramsey numbers extend this definition by assigning
a set of colors to each edge.
For integers $n\ge3$ and $r>s\ge1$, let $R(n;r,s,1)$
be the least $N$ such that every map
$\chi:E(K_N)\to\binom{[r]}{s}$ contains an $n$-vertex
clique $H$ satisfying
\(
 \bigcap_{e\in E(H)}\chi(e)\ne\varnothing.
\)
Thus, some color belongs to every edge label of $H$,
although the labels need not be identical.
In particular, $R(n;r,1,1)=R(n;r)$.
More generally, for $1\le t\le s$, let $R(n;r,s,t)$
be the least $N$ such that every such map contains
an $n$-vertex clique $H$ with
\(
 \left|\bigcap_{e\in E(H)}\chi(e)\right|\ge t.
\)

In this paper we fix the label size $s=r-1$ and study
$R(n;r,r-1,t)$ as $t$ varies from $1$ to $r-1$.
Here every edge label omits exactly one color, so $t=1$
asks for a clique on which some color is never omitted.
This endpoint was studied by Erd\H{o}s, Hajnal, and
Rado~\cite{EHR1965} in 1965.
They conjectured that $R(n;r,r-1,1)\le2^{\delta(r)n}$
for some $\delta(r)\to0$ as $r\to\infty$.
Erd\H{o}s and Szemer\'edi~\cite{ES1972} proved the conjecture
in 1972. Their upper bound, together with the random-coloring
lower bound, gives
\[
 2^{an/r}
 \le R(n;r,r-1,1)
 \le2^{bn\log r/r},
\]
where $a,b>0$ are absolute constants.

At the other endpoint, $t=r-1$ requires all edge labels
of the clique to be identical, giving
$R(n;r,r-1,r-1)=R(n;r)$.
Consequently, monotonicity and the preceding bounds yield
\[
 2^{c_1n/r}
 \le R(n;r,r-1,1)
 \le\cdots\le R(n;r,r-1,r-1)
 =R(n;r)
 \le2^{c_2nr\log r},
\]
where $c_1,c_2>0$ are absolute constants.
To describe the dependence on $t$, we write
\[
 R(n;r,r-1,t)=\exp\bigl(nf_{r,n}(t)\bigr),
 \qquad
 f(t):=f_{r,n}(t)=\frac1n\log R(n;r,r-1,t).
\]
We study the size of $f(t)$ and its variation with $t$.
Our first theorem gives bounds throughout the range
$1\le t\le r-1$.

\begin{theorem}\label{thm:stages}
There exist absolute constants $C\ge101$, $\varepsilon_3>0$,
and $0<c_0\le1$ such that, for every fixed pair of integers
$r\ge2$ and $1\le t\le r-1$, as $n\to\infty$,
\[
 \frac12\log\frac{r}{r-t}+o(1)
 \le f(t)
 \le \frac{r}{r-t}\log\frac{r}{r-t}
     +C+O_{r,t}\!\left(\frac{\log n}{n}\right).
\]
Moreover, for each fixed pair $(r,t)$, the following bounds
hold for all sufficiently large $n$.
The comparisons with $\frac1n\log R(n;k)$ hold for every
$n\ge3$.
\begin{enumerate}
\item[(i)]
For $1\le t\le r/2$,
\(
 f(t)\le
 \min\left\{
   \frac1n\log R(n;2),\,
   \frac{Ct}{r}\log\frac rt
 \right\},
\)
where $R(n;2)\le3.69507^n$ for sufficiently large $n$;
see~\cite[Theorem~1]{LuWang2026}.

\item[(ii)]
For $\frac{r}{2}<t\le 2r/3$,
$f(t)\le\frac1n\log R(n;3)$, where
\(
 R(n;3)\le
 \exp\bigl((3\log3-\varepsilon_3)n\bigr)
\)
for sufficiently large $n$;
see~\cite[Theorem~1.1]{BBCGHMST2026}.

\item[(iii)]
For $2r/3<t\le r-1$, put $h=r-t$ and write
$r=mh+a$, where $0\le a<h$.
Let $k=\lceil r/h\rceil$ and let $d=\gcd(r,h)$ be the greatest common divisor of $r$ and $h$.
Then
\[
 f(t)\le
 \min\left\{
 \begin{aligned}
 &\frac1n\log R(n;k),\\
 &\frac rh\log\frac rh
   +\min\left\{
     C,\,
     \log\frac hd-\frac ah\log\frac ah
   \right\}
   +O_{r,h}\!\left(\frac{\log n}{n}\right)
 \end{aligned}
 \right\},
\]
where
\(
 R(n;k)\le
 \exp\!\left[
   \left(k\log k-\frac{c_0}{k\log^2(2k)}\right)n
 \right]
 \;\text{for }n\ge Kk^2\log^6(2k),
\)
with $K>0$ absolute; see~\cite[Theorem~1.1]{Jo2026}.
We use the convention $0\log0=0$.
\end{enumerate}
\end{theorem}

\begin{remark}\label{rem:stagescomparison}
For $t=1$, the upper bound in
Theorem~\ref{thm:stages}(i) matches, up to an absolute constant,
the bound of Erd\H{o}s and Szemer\'edi~\cite{ES1972}:
\(
 f(1)\le C\frac{\log r}{r}.
\)
The lower bound follows from the first-moment method and gives
$f(1)\ge\frac12\log\frac{r}{r-1}+o(1)$.
Our bound extends the Erd\H{o}s--Szemer\'edi estimate
to several common colors.
For larger $t$, optimizing color partitions yields
\[
 f(t)\le
 \frac{r}{r-t}\log\frac{r}{r-t}+C+o(1),
\]
where $C$ is an absolute constant independent of $r$, $t$,
and $n$, and the error tends to zero as $n\to\infty$
with $r$ and $t$ fixed.
\end{remark}

Theorem~\ref{thm:stages} describes three regimes according to
the density $t/r$. To state them without assuming convergence
in $n$, define
\[
 \underline f_r(t):=\liminf_{n\to\infty}f_{r,n}(t),
 \qquad
 \overline f_r(t):=\limsup_{n\to\infty}f_{r,n}(t).
\]
We first let $n\to\infty$ with $r$ and $t$ fixed,
and then study how $\underline f_r(t)$ and $\overline f_r(t)$
vary with $r$ and $t$.

\begin{corollary}\label{cor:ratesbounds}
The bounds in Theorem~\ref{thm:stages} have the following consequences.
\begin{enumerate}
\item[(i)]
For all integers $r\ge2$ and $1\le t<r$,
\[
 \frac12\log\frac1{1-t/r}
 \le\underline f_r(t)\le\overline f_r(t)
 \le\frac1{1-t/r}\log\frac1{1-t/r}+C,
\]
where $C>0$ is absolute. For fixed $r$, both functions are
nondecreasing in $t$.

\item[(ii)]
Along any sequence with $r\to\infty$ and $1\le t<r$,
\[
 \begin{aligned}
 \underline f_r(t)\to0
 &\quad\Longleftrightarrow\quad t/r\to0
 \quad\Longleftrightarrow\quad \overline f_r(t)\to0,\\
 \underline f_r(t)\to\infty
 &\quad\Longleftrightarrow\quad t/r\to1
 \quad\Longleftrightarrow\quad \overline f_r(t)\to\infty.
 \end{aligned}
\]
Moreover, $c_1\le\underline f_r(t)\le\overline f_r(t)\le c_2$
eventually for some constants $0<c_1\le c_2<\infty$
if and only if $0<\liminf t/r\le\limsup t/r<1$.

\item[(iii)]
As $0<t/r\to0$, we have
$\underline f_r(t)=(t/r)^{1+o(1)}$ and
$\overline f_r(t)=(t/r)^{1+o(1)}$.
More explicitly, for every $\varepsilon>0$, there exists
$\eta>0$ such that $0<t/r<\eta$ implies
\[
 (t/r)^{1+\varepsilon}
 \le\underline f_r(t)\le\overline f_r(t)
 \le(t/r)^{1-\varepsilon}.
\]
\end{enumerate}
\end{corollary}

Corollary~\ref{cor:ratesbounds} determines when the exponential
rates $\underline f_r(t)$ and $\overline f_r(t)$
(the lower and upper limits of $f_{r,n}(t)$ as $n\to\infty$
with $r,t$ fixed) tend to zero, remain bounded above and
away from zero, or tend to infinity.
We next ask whether these rates converge as $r\to\infty$
and $\frac{t}{r}$ approaches a fixed density.
Comparisons between different numbers of colors yield
the following one-sided limits.

\begin{theorem}\label{thm:rates}
For every $\alpha\in(0,1)$, the limits
\[
 F^-(\alpha)=
 \lim_{\substack{r\to\infty,\ t/r\to\alpha\\t/r>\alpha}}
 \underline f_r(t),
 \qquad
 F^+(\alpha)=
 \lim_{\substack{r\to\infty,\ t/r\to\alpha\\t/r>\alpha}}
 \overline f_r(t)
\]
exist. The functions $F^-$ and $F^+$ are finite, positive,
nondecreasing, and right-continuous on $(0,1)$.

For every $\alpha\in(0,1)$,
$\underline f_r(t)\to F^-(\alpha)$ along every sequence
with $r\to\infty$ and $t/r\to\alpha$ if and only if
$F^-$ is continuous at $\alpha$.
The corresponding statement holds for
$\overline f_r(t)$ and $F^+$.
Consequently, each convergence statement holds for all
$\alpha\in(0,1)$ except possibly at countably many points.
\end{theorem}

\begin{remark}\label{rem:ratescomparison}
Corollary~\ref{cor:ratesbounds} controls the sizes of the rates,
whereas Theorem~\ref{thm:rates} describes their limits as the
color parameters vary. The theorem neither identifies $F^-$
with $F^+$ nor establishes convergence of $f_{r,n}(t)$
as $n\to\infty$ for fixed $r,t$.
\end{remark}

We also study the increase between consecutive values of $t$.
For $r\ge3$, define
\[
 \Delta_t:=f(t+1)-f(t)
 =\frac1n\log\frac{R(n;r,r-1,t+1)}{R(n;r,r-1,t)}
 \ge0,
 \qquad 1\le t\le r-2.
\]
The sum of these increments is $f(r-1)-f(1)$.
For fixed $r$ and $t$, if
$\liminf_{n\to\infty}\Delta_t>0$, then there exists $c>0$
such that
\[
 R(n;r,r-1,t+1)\ge 2^{cn}R(n;r,r-1,t)
\]
for all sufficiently large $n$.
We call this exponential separation.
Our third theorem gives quantitative bounds for these
increments and the contribution of the final increments.

\begin{theorem}\label{thm:separation}
Let $r,n\ge3$ and $\Delta_t=f(t+1)-f(t)$ for $1\le t\le r-2$.
\begin{enumerate}
\item[(i)] As $r\to\infty$, for every $n\ge r^2$, at least
$(1-o(1))\log r/\log\log r$ distinct indices $1\le t<r/2$
satisfy $\Delta_t\ge(1-o(1))/(2r)$.
These indices may depend on $n$.
Each $o(1)$ term depends only on $r$ and tends to zero
as $r\to\infty$.

\item[(ii)] For $r,n\ge3$ and $2\le h\le r-1$,
\[
 \frac{\sum_{t=r-h}^{r-2}\Delta_t}
      {\sum_{t=1}^{r-2}\Delta_t}
 =\frac{f(r-1)-f(r-h)}{f(r-1)-f(1)}
 \ge1-\frac1{\lfloor r/\lceil r/h\rceil\rfloor}
 \ge1-\frac1{\lceil h/2\rceil}.
\]
Thus the final $h-1$ increments account for at least
$1-1/\lceil h/2\rceil$ of the total increase.

\item[(iii)] For every fixed $r\ge4$,
\(
 \liminf_{n\to\infty}\Delta_{r-2}
 \ge\frac{\lfloor r/2\rfloor-1}{2}\log2>0.
\)

\item[(iv)] For $r,n\ge3$, if $r$ is even or $n=3$, then
$\Delta_{r-2}>\sum_{t=1}^{r-3}\Delta_t$.
\end{enumerate}
\end{theorem}

\begin{remark}\label{rem:incrementscomparison}
Part~(i) gives many separated pairs in the range $t<r/2$,
with indices that may depend on $n$.
Part~(ii) shows that the final $h-1$ increments account for
a proportion tending to one when $h\to\infty$, even if $h=o(r)$.
Parts~(iii) and~(iv) describe the final increment.
\end{remark}

The preceding results also have an equivalent formulation
in terms of cliques using few colors.
For each edge $e$, let $\psi(e)$ be the unique color omitted
from $\chi(e)$, so that $\chi(e)=[r]\setminus\{\psi(e)\}$.
Then every clique $H$ satisfies
\[
 \left|\bigcap_{e\in E(H)}\chi(e)\right|
 =r-\left|\{\psi(e):e\in E(H)\}\right|.
\]
Consequently, $R(n;r,r-1,t)$ is the least $N$ such that
every ordinary $r$-coloring of $K_N$ contains an $n$-vertex
clique using at most $r-t$ colors.

This parameter belongs to the generalized Ramsey framework
studied by Chung and Liu~\cite{ChungLiu1978}.
It is denoted by $R^r_{r-t}(K_n)$ by Budden, Stender,
and Zhang~\cite{BSZ2017}, and by $GR(r,K_n,r-t)$
by Lidick\'y et al.~\cite{LMPV2025}.
Finite bounds and exact values for these parameters appear
in~\cite{HM1999,BSZ2017,LMPV2025}.
A related inverse problem, studied by Erd\H{o}s and
Gy\'arf\'as~\cite{EG1997}, fixes the order of the complete
graph and asks for the minimum number of colors in an
edge-coloring for which every $n$-vertex clique uses
at least a fixed number of colors.

Another direction is to fix $t=1$ and vary the label size $s$,
as in the work of Conlon et al.~\cite{CFHMSV2024}
and Arag\~ao et al.~\cite{ACMMM2024}.
This leads to a different dependence on the number of colors.
For each fixed $0<\eta<1/2$, Corollary~\ref{cor:ratesbounds}
gives
\(
 R(n;r,r-1,t)=\exp(\Theta_\eta(n))
 \;\text{when }\eta\le t/r\le1-\eta,
\)
whereas the bounds of~\cite{CFHMSV2024} give
\(
 R(n;r,s,1)=\exp(\Theta_\eta(nr))
 \;\text{when }\eta\le s/r\le1-\eta.
\)
In both statements, the implicit constants depend only on
$\eta$, and $n$ is sufficiently large in terms of the
color parameters.

Section~\ref{sec:constructions} proves
Theorem~\ref{thm:stages} and Corollary~\ref{cor:ratesbounds}
using the first-moment method, sparse graph bounds, and color
partitions. Section~\ref{sec:variation} proves
Theorem~\ref{thm:rates} by comparing different numbers of colors.
Section~\ref{sec:rates} proves Theorem~\ref{thm:separation}
using disjoint intervals of indices and product colorings.
Section~\ref{sec:finite} collects elementary finite values,
including the triangle formula used in the proof of
Theorem~\ref{thm:separation}.

\section{The proof of Theorem 1}\label{sec:constructions}

The Erd\H{o}s--Szemer\'edi theorem~\cite{ES1972} states
that there exists an absolute constant $c_*>0$ such that,
for every fixed real number $k\ge2$ and all sufficiently
large $N$, every graph $H$ on $N$ vertices satisfying
\(
 e(H)\ge\left(1-\frac1k\right)\binom N2
\)
contains a clique or an independent set with at least
\(
 c_*\frac{k\log N}{\log k}
\)
vertices.

We use the following explicit version due to Liu,
which specifies the constant and the admissible range
of $k$ in terms of $N$.

\begin{lemma}[Liu~{\cite[Theorem~1]{Liu2026}}]
\label{lem:liu}
Let $H$ be a graph on $N\ge3$ vertices.
Suppose that $0<c\le1/100$ and $k$ is a real number
satisfying $2\le k\le N/(3c)$.
If
\(
 e(H)\ge\left(1-\frac1k\right)\binom N2,
\)
then $H$ contains a clique or an independent set with at least
$ck\log N/\log k$ vertices.
\end{lemma}

\begin{proof}[Proof of Theorem~\ref{thm:stages}]
We first prove the  lower bound, then the  upper bound
and the color-partition estimates.
The finite estimates will also be used in the other two proofs.

For integers $n\ge3$ and $r>h\ge1$, put
$X=\left(n!/(2\binom rh)\right)^{1/n}(r/h)^{(n-1)/2}$.

\begin{claim}\label{prop:firstmoment}
$R(n;r,r-1,r-h)\ge\lfloor X\rfloor+1>X$.
\end{claim}

\begin{claimproof}
Set $N=\lfloor X\rfloor$.
If $N<n$, then
$R(n;r,r-1,r-h)\ge n\ge N+1$, which proves the claim.
Now suppose that $N\ge n$.
Color the edges of $K_N$ independently and uniformly
from $[r]$.
For a fixed $n$-element vertex set $U$ and an
$h$-element color set $J$, all $\binom n2$ edge colors
lie in $J$ with probability
$\left(\frac{h}{r}\right)^{\binom n2}$.
Let $B$ count the pairs $(U,J)$ with this property.
Linearity of expectation gives
\[
 \mathbb E B
 =\binom Nn\binom rh\left(\frac hr\right)^{\binom n2}
 \le\frac{N^n}{n!}\binom rh\left(\frac hr\right)^{\binom n2}
 \le\frac12.
\]
The last inequality follows from $N\le X$ and the definition of $X$.
Since $B$ is a nonnegative integer,
$\mathbb P(B\ge1)\le\mathbb E B\le1/2$.
Thus a coloring with $B=0$ exists. Every clique using at most
$h$ colors would have its colors contained in an $h$-element
set $J$, so this coloring has no such $K_n$.
The omitted-color formulation gives $R(n;r,r-1,r-h)>N$.
As this Ramsey number is an integer, it is at least
$N+1=\lfloor X\rfloor+1>X$.
\end{claimproof}

For fixed $r,h$, Stirling's formula gives
\[
 n!=\sqrt{2\pi n}\left(\frac ne\right)^n
      \left(1+O\left(\frac1n\right)\right),
\]
and hence
\[
 \log(n!)
 =n\log n-n+\frac12\log(2\pi n)
  +O\left(\frac1n\right).
\]
Substituting this into the definition of $X$, we obtain
\[
 \begin{aligned}
 \frac{\log X}{n}
 &=\frac{\log(n!)-\log\left(2\binom rh\right)}{n^2}
   +\frac{n-1}{2n}\log\frac rh\\
 &=\frac12\log\frac rh
   +\frac{\log n-1-\frac12\log(r/h)}{n}
   +\frac{\frac12\log(2\pi n)-\log\left(2\binom rh\right)}{n^2}
   +O\left(\frac1{n^3}\right)\\
 &=\frac12\log\frac rh
   +\frac{\log n-1-\frac12\log(r/h)}{n}
   +O_{r,h}\left(\frac{\log n}{n^2}\right).
 \end{aligned}
\]
By Claim~\ref{prop:firstmoment},
$R(n;r,r-1,r-h)>X$, so
\[
 f_{r,n}(r-h)>\frac{\log X}{n}
 =\frac12\log\frac rh+\frac{\log n}{n}
  -O_{r,h}\left(\frac1n\right).
\]
In particular,
$f_{r,n}(r-h)\ge\frac12\log(r/h)+o(1)$ as $n\to\infty$.
Moreover, $\log n\to\infty$, whereas
$1+\frac12\log(r/h)$ is fixed and the remaining error is
$o(1/n)$. Thus $\frac{\log X}{n}>\frac12\log(r/h)$
for all sufficiently large $n$, and consequently
$f_{r,n}(r-h)>\frac12\log(r/h)$ for all such $n$.

\begin{claim}\label{prop:explicitsparse}
For all integers $n\ge3$, $r\ge2$, and $1\le t\le r/2$,
\[
 R(n;r,r-1,t)
 \le
 \max\left\{
 2n,\,
 \left\lceil\left(\frac rt\right)^{100nt/r}\right\rceil
 \right\}.
\]
\end{claim}

\begin{claimproof}
Set
\(
 N=\max\left\{
   2n,\,
   \left\lceil\left(\frac rt\right)^{100nt/r}\right\rceil
 \right\}.
\)
Then $N\ge2n$ and
\(
 \log N\ge\frac{100nt}{r}\log\frac rt.
\)

Let $\chi:E(K_N)\to\binom{[r]}{r-1}$ be an arbitrary
set-coloring. Since each edge label omits exactly one color,
there is a unique ordinary coloring $\psi:E(K_N)\to[r]$
such that $\chi(e)=[r]\setminus\{\psi(e)\}$ for every edge $e$.
For every vertex set $U$ with $|U|\ge2$,
\[
 \bigcap_{e\in E(K_N[U])}\chi(e)
 =
 [r]\setminus\{\psi(e):e\in E(K_N[U])\}.
\]
Thus it suffices to find an $n$-element vertex set whose
edges use at most $r-t$ colors under $\psi$.

For each $i\in[r]$, let
\[
 a_i=\bigl|\{e\in E(K_N):\psi(e)=i\}\bigr|.
\]
By relabeling the colors, we may assume without loss of
generality that $a_1\le a_2\le\cdots\le a_r$.
Let $G$ be the spanning subgraph of $K_N$ consisting of
the edges whose colors under $\psi$ belong to $\{1,\ldots,t\}$.
Since
\[
 \frac1t\sum_{i=1}^t a_i
 \le\frac1r\sum_{i=1}^r a_i
 =\frac1r\binom N2,
\]
we have
\[
 e(G)=\sum_{i=1}^t a_i\le\frac tr\binom N2,
 \qquad
 e(\overline G)=\binom N2-e(G)
 \ge\left(1-\frac tr\right)\binom N2.
\]
We consider two cases.

\medskip
\noindent
\textit{Case 1: $N\ge3r/(100t)$.}
Put $k=r/t$ and $c=1/100$.
Since $t\le r/2$ and $N\ge3ck$, we have
$2\le k\le N/(3c)$.
Moreover, $N\ge2n\ge6$ and
\[
 e(\overline G)\ge\left(1-\frac1k\right)\binom N2.
\]
Lemma~\ref{lem:liu} therefore gives a clique or an independent
set in $\overline G$ with at least
\[
 \frac{ck\log N}{\log k}
 =\frac{r\log N}{100t\log(r/t)}
 \ge
 \frac{r}{100t\log(r/t)}
 \cdot\frac{100nt}{r}\log\frac rt
 =n
\]
vertices.
Choose an $n$-element subset $U$ of this set.

If $U$ is a clique in $\overline G$, all edges of $K_N[U]$
have colors in $\{t+1,\ldots,r\}$ under $\psi$.
If $U$ is independent in $\overline G$, all edges of $K_N[U]$
have colors in $\{1,\ldots,t\}$.
Since $t\le r-t$, in either case these edges use at most
$r-t$ colors.

\medskip
\noindent
\textit{Case 2: $N<3r/(100t)$.}
The bound on $e(G)$ gives
\[
 e(G)\le\frac{tN(N-1)}{2r}
 <\frac{tN^2}{2r}
 <\frac{3N}{200}.
\]
For each edge of $G$, choose one endpoint, and let $S$ be
the set of chosen vertices.
Then $|S|\le e(G)$, and every edge of $G$ has an endpoint
in $S$. Hence $V(G)\setminus S$ is independent in $G$.
Since $N\ge2n$,
\[
 |V(G)\setminus S|
 =N-|S|
 \ge N-e(G)
 >\frac{197N}{200}
 >n.
\]
Choose an $n$-element subset $U$ of $V(G)\setminus S$.
All edges of $K_N[U]$ have colors in $\{t+1,\ldots,r\}$
under $\psi$, so they use at most $r-t$ colors.

\medskip
In both cases, we obtain an $n$-element vertex set $U$ with
\(
 \bigl|\{\psi(e):e\in E(K_N[U])\}\bigr|\le r-t.
\)
Returning to the original set-coloring, we have
\[
 \left|\bigcap_{e\in E(K_N[U])}\chi(e)\right|
 =
 r-\bigl|\{\psi(e):e\in E(K_N[U])\}\bigr|
 \ge t.
\]
Since $\chi$ was arbitrary, it follows that
$R(n;r,r-1,t)\le N$, as claimed.
\end{claimproof}

The following comparison gives the bounds involving ordinary
Ramsey numbers in Theorem~\ref{thm:stages}(i)--(iii).

\begin{claim}\label{cl:colorgrouping}
For integers $n\ge3$ and $1\le h<r$, let
$k=\lceil r/h\rceil$. Then
\[
 R(n;r,r-1,r-h)\le R(n;k).
\]
Consequently, for $h=r-t$,
\[
 f(t)\le
 \begin{cases}
  \dfrac1n\log R(n;2),
    & 1\le t\le r/2,\\[6pt]
  \dfrac1n\log R(n;3),
    & r/2<t\le2r/3,\\[6pt]
  \dfrac1n\log R(n;k),
    & 2r/3<t\le r-1.
 \end{cases}
\]
\end{claim}

\begin{claimproof}
Let $N=R(n;k)$, and consider an arbitrary ordinary
$r$-coloring $\psi:E(K_N)\to[r]$.
Partition $[r]$ into $k$ nonempty sets
$A_1,\ldots,A_k$, each of size at most $h$.
Define a $k$-coloring $\phi:E(K_N)\to[k]$ by
\(
 \phi(e)=j
 \quad\text{if }\psi(e)\in A_j.
\)
By the definition of $R(n;k)$, there is an $n$-vertex
clique whose edges all receive the same color $j$ under
$\phi$. All its edge colors under $\psi$ belong to $A_j$,
so it uses at most $h$ original colors.
By the omitted-color formulation, this proves
\(
 R(n;r,r-1,r-h)\le R(n;k).
\)

Taking $h=r-t$ and logarithms gives
$f(t)\le n^{-1}\log R(n;k)$.
If $1\le t\le r/2$, then $r/2\le h<r$ and $k=2$.
If $r/2<t\le2r/3$, then $r/3\le h<r/2$ and $k=3$.
In the remaining range, $h<r/3$ and $k\ge4$.
These give the stated bounds.
\end{claimproof}

We next optimize the color partition to refine the preceding bound.
Fix integers $1\le h<r$ and write $r=mh+a$, where
$m\ge1$ and $0\le a<h$.
Partition $[r]$ into nonempty groups $A_1,\ldots,A_\ell$,
and put
\[
 w_i=|A_i|,
 \qquad
 w=(w_1,\ldots,w_\ell),
 \qquad
 w_i\le h,
 \qquad
 \sum_{i=1}^{\ell}w_i=r.
\]
A group of size $h$ is called full, and a group of size one
is called a one-color group.
An index set $J\subseteq[\ell]$ is called allowed if
\[
 \left|\bigcup_{i\in J}A_i\right|
 =\sum_{i\in J}w_i\le h.
\]
The following claim describes the correspondence between
vertices and sequences of group indices.

\begin{claim}\label{cl:vertexsequences}
Let $n\ge3$, and let $\psi:E(K_N)\to[r]$ be an ordinary
coloring.
There is an injective assignment of a finite, possibly empty,
sequence $\sigma(v)$ over $[\ell]$ to each vertex $v$
with the following properties.
If $\sigma(v)=(i_1,\ldots,i_q)$, then there are distinct
vertices $v_1,\ldots,v_q,v$ such that
\[
 \sigma(v_j)=(i_1,\ldots,i_{j-1})
 \qquad(1\le j\le q),
\]
where the sequence is empty when $j=1$, and
\[
 \psi(v_jv_p)\in A_{i_j}
 \quad(1\le j<p\le q),
 \qquad
 \psi(v_jv)\in A_{i_j}
 \quad(1\le j\le q).
\]
If $\psi$ contains no $K_n$ using at most $h$ colors, then
every assigned sequence $(i_1,\ldots,i_q)$ satisfies
\[
 \bigl|\{p\in[q]:i_p\in J\}\bigr|\le n-2
 \qquad\text{for every allowed }J.
\]
\end{claim}

\begin{claimproof}
We construct the assignment recursively.
At each stage, the current vertex set carries an associated
sequence of group indices.
We start with $V(K_N)$ and the empty sequence.

Suppose that a nonempty current set $S$ has associated
sequence $\sigma=(i_1,\ldots,i_q)$, allowing $q=0$.
Choose a vertex $v\in S$ and assign $\sigma(v)=\sigma$.
Remove $v$ and partition $S\setminus\{v\}$ into the sets
\[
 S_i=\{u\in S\setminus\{v\}:\psi(vu)\in A_i\},
 \qquad i\in[\ell].
\]
Since $A_1,\ldots,A_\ell$ partition $[r]$, the sets $S_i$
are pairwise disjoint and their union is $S\setminus\{v\}$.
For each nonempty $S_i$, associate the extended sequence
\[
 \sigma i=(i_1,\ldots,i_q,i)
\]
with that set and repeat the procedure entirely within $S_i$.
Thus different parts are processed separately, and every
later subdivision of $S_i$ remains inside $S_i$.
If $S=\{v\}$, there are no nonempty parts and this branch
of the recursion terminates.

Every vertex subsequently selected within $S_i$ receives
a sequence beginning with $\sigma i$.
Only the first vertex selected from $S_i$ receives exactly
$\sigma i$; the others receive proper extensions of that
sequence.

Each selected vertex is removed and cannot be selected again.
Every unselected vertex belongs to exactly one smaller
current set at the next step along its recursive path.
Since the initial set is finite, every vertex is eventually
selected exactly once.

We next verify that the assigned sequences are distinct.
The empty sequence identifies the initial set.
Inductively, if $\sigma$ identifies at most one current set,
then $\sigma i$ identifies at most one of its parts.
Thus each sequence identifies at most one current set.
That set assigns its sequence only to its selected vertex,
while all vertices selected later within its subsets receive
longer sequences.
The assignment is therefore injective.

Now let $v$ have sequence $(i_1,\ldots,i_q)$.
For $1\le j\le q$, let $v_j$ be the vertex selected from
the set associated with $(i_1,\ldots,i_{j-1})$ on the
recursive path to $v$.
These vertices are distinct, and
\[
 \sigma(v_j)=(i_1,\ldots,i_{j-1}).
\]
At the selection of $v_j$, all later vertices
$v_{j+1},\ldots,v_q,v$ lie in the part indexed by $i_j$.
Consequently, every edge from $v_j$ to one of these vertices
has its color in $A_{i_j}$.
This proves the stated edge properties.

Finally, assume that $\psi$ contains no $K_n$ using at most
$h$ colors.
Suppose that, for some allowed set $J$, an assigned sequence
has at least $n-1$ positions $p$ satisfying $i_p\in J$.
Choose $n-1$ of the corresponding vertices $v_p$ and add $v$.
By the edge properties, every edge among these $n$ vertices
has its color in the group associated with its earlier
endpoint on the recursive path.
All these colors therefore belong to
$\bigcup_{i\in J}A_i$.
Since $J$ is allowed,
\[
 \left|\bigcup_{i\in J}A_i\right|
 =\sum_{i\in J}w_i\le h.
\]
We have obtained a $K_n$ using at most $h$ colors,
a contradiction.
Hence every assigned sequence satisfies the required
bound.
\end{claimproof}

For a sequence $\sigma=(i_1,\ldots,i_q)$, define
\[
 z_i=z_i(\sigma):=\bigl|\{p\in[q]:i_p=i\}\bigr|
 \qquad(i\in[\ell]).
\]
Thus $w_i$ is the number of colors in $A_i$, whereas $z_i$
is the number of times the index $i$ occurs along the sequence,
and $\sum_i z_i=q$.
The last conclusion of Claim~\ref{cl:vertexsequences} becomes
$\sum_{i\in J}z_i\le n-2$ for every allowed $J$.
There are exactly
$\binom{z_1+\cdots+z_\ell}{z_1,\ldots,z_\ell}$ sequences
with a prescribed occurrence vector $(z_1,\ldots,z_\ell)$.
For $n\ge3$, define
\[
 W_n(w)=
 \sum_{\substack{
   z_1,\ldots,z_\ell\in\mathbb Z_{\ge0}\\
   \sum_{i\in J}z_i\le n-2\ \text{for every allowed }J
 }}
 \binom{z_1+\cdots+z_\ell}{z_1,\ldots,z_\ell}.
\]
This counts all sequences satisfying the indicated constraints,
including the empty sequence.

To estimate $W_n(w)$, normalize the counts by
$x_i=z_i/(n-2)$ and allow arbitrary nonnegative real coordinates.
Call $x=(x_1,\ldots,x_\ell)$ feasible if
$\sum_{i\in J}x_i\le1$ for every allowed $J$.
Putting $L=\sum_i x_i$, define
\[
 \Phi(w)=
 \max_{\substack{
   x_i\ge0\\
   \sum_{i\in J}x_i\le1\ \text{for every allowed }J
 }}
 \left(L\log L-\sum_{i=1}^{\ell}x_i\log x_i\right),
 \qquad
 \Psi(r,h)=\min_w\Phi(w),
\]
where the minimum is over all positive integer group-size
vectors with $w_i\le h$ and $\sum_iw_i=r$.
We use $0\log0=0$.
Every singleton is allowed, so the feasible region is a
compact subset of $[0,1]^\ell$.
The objective is continuous, and there are finitely many
choices of $w$, so both extrema exist.

The next claim bounds the Ramsey number by counting these
sequences and then estimates the resulting coefficient
using suitable partitions.

\begin{claim}\label{thm:partitionvariational}
Let $r>h\ge1$ be integers and write $r=mh+a$, where
$m=\lfloor r/h\rfloor$ and $0\le a<h$.
For any partition $w=(w_1,\ldots,w_\ell)$ of $r$ into
positive integer group sizes at most $h$, the following hold.
\begin{enumerate}
\item[(1)] For $n\ge3$,
\[
 R(n;r,r-1,r-h)
 \le1+W_n(w)
 \le1+(n-1)^\ell\exp\bigl((n-2)\Phi(w)\bigr).
\]

\item[(2)] If $a=0$, then $\Psi(r,h)\le m\log m$.

\item[(3)] If $a>0$, then
\[
 \Psi(r,h)
 \le\frac rh\log\frac rh
   +\log\frac h{\gcd(r,h)}
   -\frac ah\log\frac ah.
\]
\end{enumerate}
\end{claim}

\begin{claimproof}
For (1), consider an ordinary $r$-coloring of $K_N$
containing no $K_n$ using at most $h$ colors.
Claim~\ref{cl:vertexsequences} assigns to each vertex a
distinct sequence counted by $W_n(w)$, so $N\le W_n(w)$.
Therefore every ordinary $r$-coloring of $K_{W_n(w)+1}$
contains a $K_n$ using at most $h$ colors.
By the omitted-color formulation,
\[
 R(n;r,r-1,r-h)\le W_n(w)+1.
\]
Since every singleton is allowed, each vector in the sum
defining $W_n(w)$ satisfies $0\le z_i\le n-2$.
There are therefore at most $(n-1)^\ell$ such vectors.
For one of them, put $q=\sum_i z_i$, $x_i=z_i/(n-2)$,
and $L=\sum_i x_i$.
If $q>0$, the multinomial theorem gives
\[
 \begin{aligned}
 \binom{q}{z_1,\ldots,z_\ell}
 &\le\prod_{i:z_i>0}\left(\frac q{z_i}\right)^{z_i}\\
 &=\exp\left(q\log q-\sum_i z_i\log z_i\right)\\
 &=\exp\left((n-2)
     \left(L\log L-\sum_i x_i\log x_i\right)\right)\\
 &\le\exp\bigl((n-2)\Phi(w)\bigr),
 \end{aligned}
\]
because $x$ is feasible.
The same bound holds when $q=0$, since that summand is one
and $\Phi(w)\ge0$.
Summing over the possible vectors proves (1).

For (2) and (3), note that the objective is nondecreasing
in each coordinate, since
\[
 \frac{\partial}{\partial x_i}
 \left(L\log L-\sum_jx_j\log x_j\right)
 =\log\frac L{x_i}\ge0
\]
for positive coordinates, with continuity handling the boundary.

For (2), take $m$ full groups.
The feasible region is $[0,1]^m$, so a maximum is attained
at $x_1=\cdots=x_m=1$.
Thus $\Phi(w)=m\log m$, proving (2).

For (3), suppose that $a>0$ and put
\(
 d=\gcd(r,h)=\gcd(a,h),\;
 b=\frac hd,\; c=\frac ad.
\)
Then $h=bd$, $a=cd$, and $1\le c<b$.
Choose the partition
\[
 w=
 \left(
 \underbrace{h,\ldots,h}_{m-1\text{ groups}},
 a,
 \underbrace{d,\ldots,d}_{b\text{ groups}}
 \right).
\]
Its entries sum to $(m-1)h+a+bd=r$.
When $m=1$, the first collection of groups is empty.

An allowed set containing a full group cannot contain any
other group. Hence the coordinate corresponding to each
full group is constrained only to lie in $[0,1]$.
By monotonicity of the objective, all these coordinates
may be set to one.

Let $u$ be the coordinate corresponding to the size-$a$
group, and let $v_1,\ldots,v_b$ be those corresponding to
the size-$d$ groups.
The feasible region is invariant under permutations of
$v_1,\ldots,v_b$.
Averaging over these permutations preserves feasibility,
because the constraints are linear, and leaves
$L=\sum_i x_i$ unchanged.
If $v=\frac1b\sum_{i=1}^b v_i$, convexity of
$x\mapsto x\log x$ gives
\(
 bv\log v\le\sum_{i=1}^b v_i\log v_i.
\)
Thus replacing $v_1,\ldots,v_b$ by their common average
cannot decrease the objective.
We may therefore assume that all size-$d$ groups have
the same coordinate $v$.

The $b$ size-$d$ groups together contain $bd=h$ colors,
so their coordinates satisfy $bv\le1$.
An allowed set containing the size-$a$ group can contain
at most $b-c$ size-$d$ groups, since
$a+(b-c)d=h$.
Conversely, these two conditions imply all constraints
involving no full group, because the coordinates are
nonnegative.
The remaining constraints are therefore exactly
\[
 0\le v\le\frac1b,\qquad
 u\ge0,\qquad
 u+(b-c)v\le1.
\]

Since the full-group coordinates contribute
$1\log1=0$, the objective becomes
\[
 L\log L-u\log u-bv\log v,
 \qquad L=m-1+u+bv.
\]
For fixed $v$, monotonicity allows us to take
\(
 u=1-(b-c)v,\; L=m+cv.
\)
Write $F(v)=L\log L-u\log u-bv\log v$ with these choices.
For $0<v\le\frac1b$, differentiation gives
\[
 \begin{aligned}
 F'(v)
 &=c(\log L+1)+(b-c)(\log u+1)-b(\log v+1)\\
 &=c\log\frac Lv+(b-c)\log\frac uv.
 \end{aligned}
\]
Here
\(
 u\ge\frac cb\ge\frac1b\ge v,
 \;  L>v,
\)
so $F'(v)>0$.
By continuity at $v=0$, the maximum is attained at
\[
 v=\frac1b,\qquad
 u=\frac cb=\frac ah,\qquad
 L=m+\frac ah=\frac rh.
\]
Consequently,
\[
 \begin{aligned}
 \Phi(w)
 &=\frac rh\log\frac rh
   -\frac ah\log\frac ah
   -b\cdot\frac1b\log\frac1b\\
 &=\frac rh\log\frac rh
   +\log\frac hd-\frac ah\log\frac ah.
 \end{aligned}
\]
Since $\Psi(r,h)\le\Phi(w)$, this proves (3).

For later use, we also note that $d\ge1$ and $a\ge1$ imply
\[
 \log\frac hd-\frac ah\log\frac ah
 =\log\frac hd+\frac ah\log\frac ha
 \le\left(1+\frac ah\right)\log h.
\]
\end{claimproof}

We next use groups whose sizes are powers of two to prove
a uniform additive bound.
Put $q=\frac rh>1$ and enlarge the absolute constant $C$
if necessary.

We begin with an auxiliary partition estimate.
For a group-size vector $w=(w_1,\ldots,w_\ell)$, let
$\Phi_\lambda(w)$ be defined in the same way as $\Phi(w)$,
except that an index set $J\subseteq[\ell]$ is allowed
when
\(
 \sum_{i\in J}w_i\le\lambda.
\)
In particular, $\Phi_h(w)=\Phi(w)$ for a partition
of the original $r$ colors.

\begin{claim}\label{cl:dyadicpartition}
Let $\lambda\ge16$ be an integer, choose an integer $d\ge0$
such that $2^d\le\lambda/16<2^{d+1}$, and put
$B_\lambda=18d+128$.
Let $M,g$ be integers with $2\lambda\le M\le4\lambda$
and $g\ge0$, and put
\[
 s=M/\lambda,
 \qquad Q=g+s.
\]
If $Q\ge\lambda^{64}$, there is a partition $w$ of
$g\lambda+M$ colors into $g$ groups of size $\lambda$
and at most $B_\lambda$ additional groups such that
\[
 \Phi_\lambda(w)\le Q\log Q+16.
\]
Each additional group has size a power of two at most
$\lambda/16$.
\end{claim}

\begin{claimproof}
Put $D=2^d$ and $v_i=2^i$ for $0\le i\le d$.
We construct $b_i$ groups of size $v_i$.
Start with $M_0=M$.
For $i<d$, choose the least integer $b_i\ge4s$ having the same
parity as $M_i$, and set $M_{i+1}=(M_i-b_i)/2$.
Finally, put $b_d=M_d$.
Thus $4s\le b_i<4s+2$ for $i<d$, and
\[
 M_i=\frac{M}{2^i}-\sum_{j<i}b_j2^{j-i}.
\]
Since $M/2^i\ge M/D\ge16s$ and
$\sum_{j<i}2^{j-i}<1$, every $M_i$ is a positive integer.
More precisely,
\[
 M_i\ge12s-2>0,
 \qquad
 \sum_{i=0}^d b_i v_i=M,
 \qquad
 12s-2\le b_d<32s.
\]
The last upper bound follows from
$b_d\le M/D=s\lambda/D<32s$.
Since $2\le s\le4$, we have $b_i<18$ for $i<d$ and $b_d<128$.
Hence the number of additional groups is at most $B_\lambda$.

Put $\mu_i=b_i/s$.
Then
\[
 \sum_{i=0}^d v_i\mu_i=\lambda,
 \qquad
 4\le\mu_i<5\quad(i<d),
 \qquad
 11\le\mu_d<32,
 \qquad
 b_i\ge2\mu_i.
\]
We next show that the allowed unions provide enough constraints
to determine the maximum of the objective.

Let
\[
 \mathcal U=
 \left\{
 u\in\mathbb Z_{\ge0}^{d+1}:
 0\le u_i\le b_i,\quad
 \sum_{i=0}^d v_i u_i=\lambda
 \right\}.
\]
Here $u_i$ counts the groups of size $v_i$ selected in a union.
Thus every vector in $\mathcal U$ describes an allowed union
whose total size is exactly $\lambda$.

We claim that every real vector $y$ satisfying
\[
 \sum_{i=0}^d v_i y_i=\lambda,
 \qquad
 |y_i-\mu_i|\le\frac12\quad(0\le i\le d)
\]
belongs to $\operatorname{conv}\mathcal U$.
We prove this by constructing a random vector
$u\in\mathcal U$ with $\mathbb Eu=y$.

If $d=0$, then $y_0=\lambda$ and $b_0=M\ge2\lambda$,
so the claim is immediate.
Suppose that $d\ge1$.
Start with the integer $T_0=\lambda$ and put $\eta_{-1}=0$.
At step $i<d$, set $c_i=y_i+\eta_{i-1}$.
We maintain the identity
\[
 T_i=c_i+\sum_{j=i+1}^d2^{j-i}y_j.
\]
Choose integers $k$ and $k+2$ of the same parity as $T_i$
such that $k\le c_i\le k+2$.
Conditioned on all earlier choices, choose
\[
 u_i=
 \begin{cases}
 k,&\text{with probability }(k+2-c_i)/2,\\
 k+2,&\text{with probability }(c_i-k)/2.
 \end{cases}
\]
These probabilities give conditional mean $c_i$ and
$|u_i-c_i|<2$ almost surely.
Now set
\[
 T_{i+1}=\frac{T_i-u_i}{2},
 \qquad
 \eta_i=\frac{c_i-u_i}{2}.
\]
The parity condition makes $T_{i+1}$ an integer.
Moreover, $|\eta_i|<1$, and the conditional mean of $\eta_i$
is zero.
The identity for the next step is
\[
 T_{i+1}
 =y_{i+1}+\eta_i+
   \sum_{j=i+2}^d2^{j-i-1}y_j,
\]
so the construction can continue.

Finally, set $u_d=T_d=y_d+\eta_{d-1}$.
The recurrence $T_i=u_i+2T_{i+1}$ gives
$\sum_i2^i u_i=\lambda$.
Since every $\eta_i$ has expectation zero, we also obtain
$\mathbb Eu_i=y_i$ for every $i$.
Furthermore,
\[
 |u_i-y_i|<3\quad(i<d),
 \qquad
 |u_d-y_d|<1.
\]
Together with $|y_i-\mu_i|\le1/2$, $\mu_i\ge4$, and
$b_i\ge2\mu_i$, this implies
\[
 0<\mu_i-\frac72<u_i
 <\mu_i+\frac72<2\mu_i\le b_i.
\]
Thus every outcome of positive probability belongs to
$\mathcal U$.
There are only finitely many outcomes, so their expectation
$y$ is a convex combination of vectors in $\mathcal U$.
This proves the claim about $\operatorname{conv}\mathcal U$.

We now compute the coefficient of the auxiliary partition.
A group of size $\lambda$ can occur only alone in an allowed union.
By coordinatewise monotonicity, its coordinate may therefore
be set to one at a maximum.
Averaging the coordinates of equal-sized groups preserves
feasibility and their total coordinate sum.
By convexity of $x\log x$, this averaging cannot decrease
the objective.
We may therefore assign a common coordinate $x_i$ to the
$b_i$ groups of size $v_i$.

The resulting objective is
\[
 F(x)=
 \left(g+\sum_{i=0}^d b_i x_i\right)
 \log\left(g+\sum_{i=0}^d b_i x_i\right)
 -\sum_{i=0}^d b_i x_i\log x_i.
\]
For positive coordinates, put $L=g+\sum_i b_i x_i$.
For every real vector $\xi=(\xi_0,\ldots,\xi_d)$,
Cauchy--Schwarz gives
\[
 \xi^{\mathsf T}\nabla^2F(x)\xi
 =\frac{(\sum_i b_i\xi_i)^2}{L}
  -\sum_i\frac{b_i\xi_i^2}{x_i}
 \le0.
\]
Thus $F$ is concave, with continuity extending this property
to the boundary.

Consider $p_i=v_i/\lambda$.
This vector is feasible because, for any allowed union of
additional groups, its coordinate sum equals its total group
size divided by $\lambda$.
Also, $\sum_i b_i p_i=s$, and hence
\[
 F(p)=Q\log Q+E,
 \qquad
 E=\sum_{i=0}^d b_i p_i\log\frac1{p_i}.
\]

To estimate $E$, put $\rho=D/\lambda$.
Then $1/32<\rho\le1/16$ and
$p_{d-j}=\rho2^{-j}$.
Using $b_i<4s+2$ for $i<d$ and
$\sum_{j\ge1}j2^{-j}=2$, we obtain
\[
 \begin{aligned}
 E
 &=s\log\frac1\rho
   +\rho\log2\sum_{j=1}^d b_{d-j}j2^{-j}\\
 &\le
   \left(5s+\frac{4s+2}{8}\right)\log2\\
 &=\left(\frac{11s}{2}+\frac14\right)\log2
 \le\frac{89}{4}\log2
 <16.
 \end{aligned}
\]

A bound on $F(p)$ alone does not bound the maximum.
We therefore verify that $p$ is a maximizer.
Its partial derivatives are
\[
 \frac{\partial F}{\partial x_i}(p)
 =b_i\left(\log Q+\log\frac1{p_i}\right).
\]
Define
\[
 y_i=
 \frac{b_i\bigl(\log Q+\log(1/p_i)\bigr)}
      {s\log Q+E}.
\]
Since $v_i=\lambda p_i$, we have
$\sum_i v_i y_i=\lambda$.
Moreover,
\[
 y_i-\mu_i
 =\mu_i\,
   \frac{\log(1/p_i)-E/s}{\log Q+E/s}.
\]
The numbers $\log(1/p_i)$ lie in $[0,\log\lambda]$.
Their weighted average, with weights $b_i p_i/s$, is $E/s$.
Thus $E/s$ lies in the same interval.
Since $\mu_i<32$ and $Q\ge\lambda^{64}$,
\[
 |y_i-\mu_i|
 \le\frac{32\log\lambda}{\log Q}
 \le\frac12.
\]
The convex-hull property proved above gives
$y\in\operatorname{conv}\mathcal U$.

For every feasible $x$ and every $u\in\mathcal U$, the
corresponding union is allowed, so $u\cdot x\le1$.
On the other hand, $u\cdot p=1$.
Taking convex combinations yields $y\cdot(x-p)\le0$.
Since $\nabla F(p)=(s\log Q+E)y$, concavity gives
\[
 F(x)
 \le F(p)+\nabla F(p)\cdot(x-p)
 \le F(p)
 <Q\log Q+16.
\]
The same inequality holds at boundary points by continuity,
with the convention $0\log0=0$.
Consequently, $\Phi_\lambda(w)\le Q\log Q+16$.
\end{claimproof}

We now apply Claim~\ref{cl:dyadicpartition} to obtain
an absolute additive bound for $\Psi(r,h)$.

\begin{claim}\label{thm:uniformconstant}
There exists an absolute constant $C\ge101$ such that,
for all integers $1\le h<r$, writing $q=r/h$, we have
\[
 \Psi(r,h)\le q\log q+C.
\]
\end{claim}

\begin{claimproof}
Put $q=r/h$ and write $r=mh+a$, where
$m=\lfloor r/h\rfloor$ and $0\le a<h$.
If $a=0$, Claim~\ref{thm:partitionvariational}(2) gives
$\Psi(r,h)\le m\log m=q\log q$.
We may therefore assume that $a>0$.

Choose an absolute constant $Q_0$ sufficiently large.
For $q\ge Q_0$, put $\lambda=\lceil(\log q)^2\rceil$.
Since $B_\lambda=O(\log\lambda)$, we may choose $Q_0$ so that
\[
 \lambda\ge16,
 \qquad
 B_\lambda+1\le\lambda,
 \qquad
 m\ge2,
 \qquad
 q\ge\lambda^{64}.
\]
We consider three ranges for $h$.

If $h<16$, then $a>0$ implies $2\le h\le15$.
Claim~\ref{thm:partitionvariational}(3) gives
\[
 \begin{aligned}
 \Psi(r,h)
 &\le q\log q+
   \log\frac{h}{\gcd(r,h)}
   -\frac ah\log\frac ah\\
 &\le q\log q+\left(1+\frac ah\right)\log h\\
 &\le q\log q+2\log15
 <q\log q+16.
 \end{aligned}
\]
Here we used $\gcd(r,h)\ge1$ and $a\ge1$.

If $16\le h\le\lambda$, apply
Claim~\ref{cl:dyadicpartition} with allowed total size $h$,
with $M=2h+a$ and $g=m-2$.
Then $2h\le M<3h$ and
\[
 Q=g+M/h=q\ge\lambda^{64}\ge h^{64}.
\]
The resulting partition has total size $r$ and gives
$\Psi(r,h)\le q\log q+16$.

It remains to consider $h>\lambda$.
Put
\[
 M=\left\lceil\left(2+\frac ah\right)\lambda\right\rceil
   +B_\lambda,
 \qquad
 s=\frac M\lambda,
 \qquad
 \delta=s-2-\frac ah.
\]
Then
\[
 2\lambda\le M\le4\lambda,
 \qquad
 0\le\delta\le\frac{B_\lambda+1}{\lambda}\le1.
\]
Apply Claim~\ref{cl:dyadicpartition} with $g=m-2$.
Its parameter is $Q=g+s=q+\delta\ge\lambda^{64}$.
It gives $g$ full groups and $b\le B_\lambda$ additional
groups whose total size is $M$.

We turn this auxiliary partition into a partition of the
original $r$ colors.
For each additional group of auxiliary size $v$, choose an
actual integer size between $1$ and $\lfloor hv/\lambda\rfloor$.
These intervals are nonempty because $h>\lambda$.
The sum of their upper endpoints satisfies
\[
 \begin{aligned}
 \sum_{\text{additional groups}}
       \left\lfloor\frac{hv}{\lambda}\right\rfloor
 &\ge\frac h\lambda M-b\\
 &\ge2h+a+\frac h\lambda B_\lambda-b\\
 &\ge2h+a.
 \end{aligned}
\]
The sum of their lower endpoints is
$b\le B_\lambda<\lambda<h<2h+a$.
Starting with size one in every group and increasing the sizes
one at a time, we can therefore make the total exactly $2h+a$.
Give each of the $m-2$ full groups actual size $h$.
The resulting partition has total size $r$.

Every actual group size is at most $h/\lambda$ times its
auxiliary size, including the full groups.
Thus a union with auxiliary total size at most $\lambda$
has actual total size at most $h$.
Every auxiliary allowed union is therefore also an actual
allowed union.
Consequently, the actual partition satisfies all the auxiliary
constraints and possibly further constraints.
Its feasible region is contained in the auxiliary feasible
region, while the objective is the same function of the group
coordinates.
Its coefficient is therefore no larger, giving
\[
 \Psi(r,h)\le(q+\delta)\log(q+\delta)+16.
\]

The change from $q$ to $q+\delta$ has bounded cost.
Indeed,
\[
 \begin{aligned}
 0
 &\le(q+\delta)\log(q+\delta)-q\log q\\
 &=\int_q^{q+\delta}(1+\log u)\,du\\
 &\le\frac{B_\lambda+1}{\lambda}\bigl(1+\log(q+1)\bigr)\\
 &=O\left(\frac{\log\log q}{\log q}\right).
 \end{aligned}
\]
The implied constant is absolute, since
$\lambda=\lceil(\log q)^2\rceil$ and
$B_\lambda=O(\log\lambda)$.
After increasing $Q_0$ if necessary, this difference is at
most one for every $q\ge Q_0$.
Combining the three ranges gives
$\Psi(r,h)\le q\log q+17$ whenever $q\ge Q_0$.

Finally, suppose that $1<q<Q_0$.
Use $m$ groups of size $h$ and one group of size $a$.
There are $k=\lceil q\rceil$ groups.
Every feasible coordinate is at most one, so coordinatewise
monotonicity bounds the coefficient by $k\log k$.
Since $k\le\lceil Q_0\rceil$, choosing
\[
 C=\max\left\{
 101,\,
 \lceil Q_0\rceil\log\lceil Q_0\rceil
 \right\}
\]
gives $\Psi(r,h)\le q\log q+C$ in this range as well.
The constant $C$ is independent of $r$ and $h$, as required.
\end{claimproof}

We now complete the proof of Theorem~\ref{thm:stages}.
Fix $r$ and $t$, and put
\[
 h=r-t,\qquad k=\lceil r/h\rceil.
\]
Let $C\ge101$ be the absolute constant in
Claim~\ref{thm:uniformconstant}.
By Claim~\ref{prop:firstmoment} and Stirling's formula,
\[
 f(t)\ge\frac12\log\frac rh+o(1)
 \qquad\text{as }n\to\infty.
\]

For the upper bound, choose a partition $w$ attaining
$\Psi(r,h)$, and let $\ell\le r$ be its number of groups.
Since $\Psi(r,h)\ge0$, Claim~\ref{thm:partitionvariational}(1)
gives
\[
 \begin{aligned}
 f(t)
 &\le\frac1n\log\left(
      1+(n-1)^\ell\exp\bigl((n-2)\Psi(r,h)\bigr)
    \right)\\
 &\le\Psi(r,h)+\frac{r\log(n-1)+\log2}{n}.
 \end{aligned}
\]
Together with Claim~\ref{thm:uniformconstant}, this yields
\[
 f(t)\le\frac rh\log\frac rh
       +C+O_{r,h}\!\left(\frac{\log n}{n}\right),
\]
which is the upper bound at the start of the theorem.
We now prove the three additional bounds.

For (i), $1\le t\le r/2$ gives $k=2$, so
Claim~\ref{cl:colorgrouping} yields
$f(t)\le n^{-1}\log R(n;2)$.
Put $\beta=(t/r)\log(r/t)>0$.
For all sufficiently large $n$, depending on $r,t$,
we have $e^{100n\beta}\ge2n$ and $\log2/n\le\beta$.
Claim~\ref{prop:explicitsparse} therefore gives
\[
 R(n;r,r-1,t)
 \le\left\lceil e^{100n\beta}\right\rceil
 \le2e^{100n\beta},
\]
and hence
\[
 f(t)\le100\beta+\frac{\log2}{n}
       \le101\beta
       \le\frac{Ct}{r}\log\frac rt.
\]
Combining the two upper bounds proves (i).
The stated estimate $R(n;2)\le3.69507^n$ for sufficiently
large $n$ is due to Lu and Wang~\cite{LuWang2026}.

For (ii), $r/2<t\le2r/3$ gives $r/3\le h<r/2$
and $k=3$.
Thus Claim~\ref{cl:colorgrouping} gives
$f(t)\le n^{-1}\log R(n;3)$.
The estimate
\[
 R(n;3)\le\exp\bigl((3\log3-\varepsilon_3)n\bigr)
\]
for an absolute constant $\varepsilon_3>0$ and all
sufficiently large $n$ follows from~\cite{BBCGHMST2026}.
This proves (ii).

For (iii), $2r/3<t\le r-1$ gives $h<r/3$ and $k\ge4$.
Claim~\ref{cl:colorgrouping} again gives
$f(t)\le n^{-1}\log R(n;k)$.
By~\cite{Jo2026}, there are absolute constants
$0<c_0\le1$ and $K>0$ such that
\[
 R(n;k)\le
 \exp\!\left[
   \left(k\log k-\frac{c_0}{k\log^2(2k)}\right)n
 \right]
 \qquad\text{for }n\ge Kk^2\log^6(2k).
\]
Since $r,h$, and hence $k$, are fixed, this condition holds
for all sufficiently large $n$.

For the other upper bound in (iii), write $r=mh+a$,
where $0\le a<h$, and put $d=\gcd(r,h)$.
If $a>0$, Claim~\ref{thm:partitionvariational}(3) gives
\[
 \Psi(r,h)\le\frac rh\log\frac rh
       +\log\frac hd-\frac ah\log\frac ah.
\]
The same inequality holds when $a=0$ by
Claim~\ref{thm:partitionvariational}(2), since then $d=h$
and the correction term is zero under the convention
$0\log0=0$.
Combining this with Claim~\ref{thm:uniformconstant} and
the partition estimate established above, we obtain
\[
 f(t)\le\frac rh\log\frac rh
   +\min\left\{
      C,\,
      \log\frac hd-\frac ah\log\frac ah
    \right\}
   +O_{r,h}\!\left(\frac{\log n}{n}\right).
\]
Taking the minimum with $n^{-1}\log R(n;k)$ proves (iii)
and completes the proof.
\end{proof}

\begin{remark}
\label{rem:ordinaryendpoints}
The endpoint identity gives
$R(n;2,1,1)=R(n;2)$ and $R(n;3,2,2)=R(n;3)$.
For $r=2,3$, the first-moment and partition bounds yield
\[
 \left(\frac12\log r+o(1)\right)n
 \le\log R(n;r)
 \le\left(r\log r+o(1)\right)n.
\]
These match the classical exponential coefficients in
the random-coloring lower bound and the
Erd\H{o}s--Szekeres upper bound.

The cited results on ordinary Ramsey numbers give the
following stronger estimates for sufficiently large $n$:
\[
 \begin{aligned}
 \left(\frac12\log2+o(1)\right)n
 &\le\log R(n;2)\le n\log3.69507,\\
 n\log2-O(1)
 &\le\log R(n;3)\le(3\log3-\varepsilon_3)n.
 \end{aligned}
\]
The two-color upper bound is due to
Lu and Wang~\cite[Theorem~1]{LuWang2026}.
For three colors, the lower bound follows from
Brada\v{c}~\cite[Theorem~1.6]{Bradac2026}, and the upper
bound follows from Balister et al.~\cite{BBCGHMST2026}.
Alternatively, Jo's bound~\cite[Theorem~1.1]{Jo2026}
permits the choice
$\varepsilon_3=\frac{c_0}{3\log^2 6}$.

These endpoint improvements are supplied by the cited results.
Our focus is the behavior at intermediate common-color
thresholds and the differences between consecutive values
of $f(t)$.
\end{remark}

We now deduce the statements about the lower and upper limits in $n$.

\begin{proof}[Proof of Corollary~\ref{cor:ratesbounds}]
For fixed $r,t$, put $\alpha=t/r$.
Taking the lower and upper limits as $n\to\infty$ in
Theorem~\ref{thm:stages} gives
\[
 \frac12\log\frac1{1-\alpha}
 \le\underline f_r(t)\le\overline f_r(t)
 \le\frac1{1-\alpha}\log\frac1{1-\alpha}+C.
\]
Since $R(n;r,r-1,t)$ is nondecreasing in $t$, both
$\underline f_r(t)$ and $\overline f_r(t)$ are nondecreasing
in $t$ for fixed $r$. This proves (i).

For (ii), consider a sequence with $r\to\infty$ and
$1\le t<r$.
When $0<\alpha\le1/2$, Theorem~\ref{thm:stages}(i)
and the inequality $-\log(1-\alpha)\ge\alpha$ give
\[
 \frac{\alpha}{2}
 \le\underline f_r(t)\le\overline f_r(t)
 \le C\alpha\log(1/\alpha).
\]
Thus $\alpha\to0$ implies that both quantities tend to zero.
Conversely, if $\alpha$ does not tend to zero, then
$\alpha\ge\eta$ along a subsequence for some $\eta>0$.
The lower bound in (i) then bounds both quantities away
from zero on that subsequence.

If $\alpha\to1$, the lower bound in (i) tends to infinity,
so both quantities tend to infinity.
Conversely, if $\alpha$ does not tend to one, then
$\alpha\le1-\eta$ along a subsequence for some $\eta>0$.
The upper bound in (i) bounds both quantities above on
that subsequence.

If $\eta\le\alpha\le1-\eta$ eventually for some $\eta>0$,
part (i) bounds both quantities above and away from zero.
Conversely, if $\liminf\alpha=0$ or $\limsup\alpha=1$,
there is a subsequence along which $\alpha$ tends to zero
or one, respectively.
The preceding arguments show that both quantities then
tend to zero or infinity along that subsequence.
This proves all assertions in (ii).

For (iii), fix $\varepsilon>0$.
Choose $\eta\in(0,1/2)$ sufficiently small that
\[
 \alpha^\varepsilon\le\frac12,
 \qquad
 C\alpha^\varepsilon\log(1/\alpha)\le1
 \qquad\text{whenever }0<\alpha<\eta.
\]
The bounds established above then give
\[
 \alpha^{1+\varepsilon}
 \le\frac{\alpha}{2}
 \le\underline f_r(t)\le\overline f_r(t)
 \le C\alpha\log(1/\alpha)
 \le\alpha^{1-\varepsilon}.
\]
This proves (iii).
\end{proof}

\section{The proof of Theorem 2}\label{sec:variation}

\begin{proof}[Proof of Theorem~\ref{thm:rates}]
Corollary~\ref{cor:ratesbounds} gives the bounds and
monotonicity of $\underline f_r(t)$ and $\overline f_r(t)$.
To establish the density limits, we also need the following
comparisons between different numbers of colors.

\begin{claim}\label{cl:palette}
For integers $n\ge3$, $r\ge2$, $1\le t<r$, and $k\ge1$,
\[
 \begin{aligned}
 &R(n;r+1,r,t)
 \le R(n;r,r-1,t)
 \le R(n;r+1,r,t+1),\\
 &R(n;kr,kr-1,kt)
 \le R(n;r,r-1,t).
 \end{aligned}
\]
\end{claim}

\begin{claimproof}
We use the omitted-color formulation throughout.

For the first inequality, merge two colors in an ordinary
$(r+1)$-coloring to obtain an $r$-coloring.
A clique using at most $r-t$ colors after merging uses
at most $r-t+1=(r+1)-t$ colors before merging.
Hence
\[
 R(n;r+1,r,t)\le R(n;r,r-1,t).
\]

For the second inequality, view an ordinary $r$-coloring
as an $(r+1)$-coloring with one unused color.
The permitted number of colors is $r-t$ in both problems,
since $(r+1)-(t+1)=r-t$.
Therefore
\[
 R(n;r,r-1,t)\le R(n;r+1,r,t+1).
\]

Finally, partition the colors of an ordinary $kr$-coloring
into $r$ groups of size $k$.
The resulting $r$-coloring contains a clique using at most
$r-t$ groups whenever the number of vertices is at least
$R(n;r,r-1,t)$.
Such a clique uses at most $k(r-t)=kr-kt$ original colors,
giving
\[
 R(n;kr,kr-1,kt)\le R(n;r,r-1,t).
\]
\end{claimproof}

For the remainder of the proof, let $g_r(t)$ denote either
$\underline f_r(t)$ or $\overline f_r(t)$, with the same
choice throughout.
Applying $n^{-1}\log$ to the inequalities in
Claim~\ref{cl:palette} and taking the corresponding lower
or upper limits as $n\to\infty$, we obtain
\[
 g_{r+1}(t)\le g_r(t)\le g_{r+1}(t+1),
 \qquad
 g_{kr}(kt)\le g_r(t).
\]
For $0<\alpha<1$, define
\[
 F(\alpha)=
 \inf\left\{
   g_q(p):q,p\in\mathbb Z,\ q\ge2,\ 1\le p<q,\ p/q>\alpha
 \right\}.
\]
The defining set is nonempty, and Corollary~\ref{cor:ratesbounds}
gives
\[
 0<\frac12\log\frac1{1-\alpha}\le F(\alpha)<\infty.
\]
Moreover, $F$ is nondecreasing by definition.
Thus its left limit
$F(\alpha-):=\lim_{\beta\uparrow\alpha}F(\beta)$ exists
for every $\alpha\in(0,1)$.

\begin{claim}\label{cl:densitylimits}
For every $\alpha\in(0,1)$ and every sequence with
$r\to\infty$, $1\le t<r$, and $t/r\to\alpha$,
\[
 F(\alpha-)
 \le\liminf g_r(t)
 \le\limsup g_r(t)
 \le F(\alpha).
\]
If $t/r>\alpha$ eventually, then $g_r(t)\to F(\alpha)$.
The function $F$ is right-continuous on $(0,1)$,
and each endpoint in the display is realized as a limit
along a suitable sequence with $r\to\infty$ and $t/r\to\alpha$.
\end{claim}

\begin{claimproof}
Fix $\alpha\in(0,1)$ and a sequence as in the statement.
Choose integers $q\ge2$ and $1\le p<q$ with $p/q>\alpha$,
and put $k=\lfloor r/q\rfloor$.
Since $kp/r\to p/q>\alpha$, we have $t\le kp$ for all
sufficiently large $r$.
Also, $k\ge1$ and $kp<kq\le r$.
Monotonicity in $t$ and the preceding comparisons give
\[
 g_r(t)\le g_r(kp)\le g_{kq}(kp)\le g_q(p).
\]
Taking the upper limit and then the infimum over $q,p$
proves $\limsup g_r(t)\le F(\alpha)$.

For any fixed $\beta\in(0,\alpha)$, we eventually have
$t/r>\beta$, so $g_r(t)\ge F(\beta)$ by definition.
Taking the lower limit and then letting $\beta\uparrow\alpha$
gives $\liminf g_r(t)\ge F(\alpha-)$.
If $t/r>\alpha$ eventually, the stronger bound
$g_r(t)\ge F(\alpha)$ gives convergence to $F(\alpha)$.
In particular, taking $t=\lfloor\alpha r\rfloor+1$
for sufficiently large $r$ realizes the upper endpoint.

To prove right-continuity, fix $\varepsilon>0$ and choose
$q,p$ with $p/q>\alpha$ and
$g_q(p)<F(\alpha)+\varepsilon$.
For every $\alpha<\beta<p/q$,
\[
 F(\alpha)\le F(\beta)\le g_q(p)<F(\alpha)+\varepsilon.
\]
Since $\varepsilon$ is arbitrary, $F$ is right-continuous
at $\alpha$.

It remains to realize the lower endpoint.
If $F(\alpha-)=F(\alpha)$, the sequence already constructed
realizes both endpoints.
Otherwise, put $\delta=F(\alpha)-F(\alpha-)>0$.
For all sufficiently large integers $j$, we have
$\alpha-1/j>0$ and $1/j<\delta$.
Choose integers $q_j\ge2$ and $1\le p_j<q_j$ such that
\[
 \frac{p_j}{q_j}>\alpha-\frac1j,
 \qquad
 g_{q_j}(p_j)
 <F\!\left(\alpha-\frac1j\right)+\frac1j
 <F(\alpha).
\]
The last inequality forces $p_j/q_j\le\alpha$;
otherwise the definition of $F(\alpha)$ would give
$g_{q_j}(p_j)\ge F(\alpha)$.
Put $r_j=jq_j$ and $t_j=jp_j$.
Then $r_j\to\infty$ and
\[
 \alpha-\frac1j<\frac{t_j}{r_j}\le\alpha,
\]
so $t_j/r_j\to\alpha$.
By the definition of $F$ and the scaling comparison,
\[
 F\!\left(\alpha-\frac1j\right)
 \le g_{r_j}(t_j)
 \le g_{q_j}(p_j)
 <F\!\left(\alpha-\frac1j\right)+\frac1j.
\]
Hence $g_{r_j}(t_j)\to F(\alpha-)$, as required.
\end{claimproof}

Apply Claim~\ref{cl:densitylimits} first with
$g_r(t)=\underline f_r(t)$ and then with
$g_r(t)=\overline f_r(t)$.
The one-sided convergence identifies the resulting functions
with $F^-$ and $F^+$, respectively.
Both functions are finite, positive, nondecreasing,
and right-continuous.

Since $F$ is right-continuous, it is continuous at $\alpha$
if and only if $F(\alpha-)=F(\alpha)$.
At such a point, the bounds in Claim~\ref{cl:densitylimits}
give convergence to $F(\alpha)$ along every sequence with
$r\to\infty$ and $t/r\to\alpha$.
At a discontinuity, the two sequences constructed in that
claim have different limits, so this unrestricted convergence
fails.

Finally, the nonempty intervals $(F(\alpha-),F(\alpha))$
associated with distinct discontinuities are pairwise disjoint,
because $F$ is nondecreasing.
Choosing a rational number in each interval shows that
there are at most countably many discontinuities.
Applying this observation to $F^-$ and $F^+$ completes
the proof.
\end{proof}

\begin{remark}\label{rem:orderlimits}
The limits as $n\to\infty$ with $r,t$ fixed cannot, in
general, be replaced by simultaneous limits in $n,r,t$.
For example, for each $n\ge3$, take
\[
 r=2\binom n2,
 \qquad
 t=\binom n2.
\]
Every ordinary $r$-coloring of $K_n$ uses at most
$\binom n2=r-t$ colors.
By the omitted-color formulation, $R(n;r,r-1,t)=n$.
Consequently,
\[
 \frac tr=\frac12,
 \qquad
 f_{r,n}(t)=\frac{\log n}{n}\longrightarrow0.
\]
This is consistent with Corollary~\ref{cor:ratesbounds}
and Theorem~\ref{thm:rates}, where the lower and upper
limits in $n$ are taken with $r,t$ fixed before the
color parameters are allowed to vary.
\end{remark}

\section{The proof of Theorem 3}\label{sec:rates}

\begin{proof}[Proof of Theorem~\ref{thm:separation}]
We first locate large increments in disjoint intervals of indices.
We then use product colorings to control the increase near the
end of the sequence.

For the next claim, fix an integer $r\ge10000$ and put
\[
 b=\left\lceil100(\log r)(\log\log r)\right\rceil,
 \qquad
 c=\frac12-\frac1{\log\log r},
 \qquad
 J=\left\lfloor\frac{\log(r/2)}{\log b}\right\rfloor.
\]

\begin{claim}\label{prop:manyinteriorincrements}
We have $J\ge1$.
For every integer $n\ge r^2$, each interval
$[b^j,b^{j+1}-1]$, $0\le j<J$, contains an index $t$
such that $\Delta_t>c/r$.
The chosen indices may depend on $n$.
Moreover, as $r\to\infty$,
\[
 J=(1+o(1))\frac{\log r}{\log\log r}.
\]
\end{claim}

\begin{claimproof}
We first verify that $2\le b\le r/2$ and $c>1/22$.
The function $x/2-100(\log x)(\log\log x)-1$ is positive
at $x=10000$ and has positive derivative
$1/2-100(1+\log\log x)/x$ for $x\ge10000$.
This gives the required bound on $b$.
Also, $\log\log r\ge\log\log10000>11/5$, so $c>1/22$.
Consequently, $J\ge1$, and the intervals
$[b^j,b^{j+1}-1]$, $0\le j<J$, are nonempty,
pairwise disjoint, and contained in $[1,r/2)$.

Fix $n\ge r^2$ and $0\le j<J$, and put
$v=b^j$ and $u=bv$.
Then $1\le v<u\le r/2$.
The exponential term in Claim~\ref{prop:explicitsparse}
exceeds $2n$, because
\[
 \frac{100nv}{r}\log\frac rv
 \ge\frac{100n\log2}{r}
 >\log(2n).
\]
The last inequality holds at $n=r^2$, and the difference
between its two sides is increasing in $n\ge r^2$.
Using $\lceil x\rceil\le2x$ for $x\ge1$ and the definitions
of $b$ and $c$, we obtain
\[
 \begin{aligned}
 f(v)
 &\le\frac{100v}{r}\log\frac rv+\frac{\log2}{n}\\
 &\le\left(\frac12-c\right)\frac ur+\frac{\log2}{n}.
 \end{aligned}
\]

At the other endpoint, Claim~\ref{prop:firstmoment} gives
\[
 \begin{aligned}
 f(u)
 &>\frac{n-1}{2n}\log\frac r{r-u}
   +\frac{\log(n!)-\log\bigl(2\binom ru\bigr)}{n^2}\\
 &\ge\frac{u}{2r}-\frac{\log2}{2n}
                    -\frac{(r+1)\log2}{n^2}.
 \end{aligned}
\]
Here we used $u/r\le1/2$,
$u/r\le-\log(1-u/r)\le\log2$,
$\log(n!)\ge0$, and $\binom ru\le2^r$.
The total error satisfies
\[
 \begin{aligned}
 E&:=\frac{3\log2}{2n}+\frac{(r+1)\log2}{n^2}\\
 &<\frac{3}{2r^2}+\frac2{r^3}
 <\frac1{22r}
 <\frac cr
 \le\frac{cv}{r}.
 \end{aligned}
\]
Therefore
\[
 f(u)-f(v)>\frac{cu}{r}-E>\frac{c(u-v)}{r}.
\]
Since
\[
 f(u)-f(v)=\sum_{t=v}^{u-1}\Delta_t,
\]
at least one of these $u-v$ increments exceeds $c/r$.
This proves the finite assertion.

Finally,
\[
 \log b=\log\log r+O(\log\log\log r+1),
\]
and hence
\[
 J=\frac{\log(r/2)}{\log b}+O(1)
  =(1+o(1))\frac{\log r}{\log\log r}.
\]
\end{claimproof}

The preceding claim proves part~(i), since
$c/r=(1-o(1))/(2r)$ and
$J=(1+o(1))\log r/\log\log r$.
The estimates hold uniformly for $n\ge r^2$, with all
asymptotic errors depending only on $r$.
We next establish the product comparisons needed for
the remaining parts.

\begin{claim}\label{prop:terminalproduct}
For integers $r,n\ge3$ and $2\le h\le r-1$, put
$k=\lceil r/h\rceil$ and $q=\lfloor r/k\rfloor$.
Then
\[
 \begin{aligned}
 &R(n;r)-1
 \ge\bigl(R(n;r,r-1,r-h)-1\bigr)^q,\\
 &R(n;r)-1
 \ge\bigl(R(n;r,r-1,r-2)-1\bigr)
     \bigl(R(n;\lfloor r/2\rfloor)-1\bigr).
 \end{aligned}
\]
\end{claim}

\begin{claimproof}
The ordinary Ramsey product inequality states that
\[
 R(n;a+b)-1
 \ge\bigl(R(n;a)-1\bigr)\bigl(R(n;b)-1\bigr)
 \qquad(a,b\ge1);
\]
see~\cite{Abbott1972}.
To see this, take an $a$-coloring on a vertex set $U$
and a $b$-coloring on a vertex set $V$, where
$|U|=R(n;a)-1$ and $|V|=R(n;b)-1$.
Choose both colorings to contain no monochromatic $K_n$
and to use disjoint sets of colors.

Color the edges on $U\times V$ using the first coloring
when the first coordinates differ, and the second coloring
when they agree.
A monochromatic clique in the first set of colors has
distinct first coordinates and projects to a monochromatic
clique in $U$.
A monochromatic clique in the second set of colors has
a common first coordinate and projects to a monochromatic
clique in $V$.
Thus neither type can have $n$ vertices, proving the product
inequality.
The construction also applies to a one-color factor,
using $R(n;1)=n$.

By Claim~\ref{cl:colorgrouping},
\[
 R(n;r,r-1,r-h)\le R(n;k).
\]
Since $qk\le r$, repeated application of the product inequality,
with the remaining colors unused, gives
\[
 \begin{aligned}
 R(n;r)-1
 &\ge\bigl(R(n;k)-1\bigr)^q\\
 &\ge\bigl(R(n;r,r-1,r-h)-1\bigr)^q.
 \end{aligned}
\]
This proves the first comparison.

For the second, put
$a=\lceil r/2\rceil$ and $b=\lfloor r/2\rfloor$.
Taking $h=2$ in Claim~\ref{cl:colorgrouping} gives
$R(n;r,r-1,r-2)\le R(n;a)$.
Since $a+b=r$, the product inequality yields
\[
 \begin{aligned}
 R(n;r)-1
 &\ge\bigl(R(n;a)-1\bigr)\bigl(R(n;b)-1\bigr)\\
 &\ge\bigl(R(n;r,r-1,r-2)-1\bigr)
      \bigl(R(n;\lfloor r/2\rfloor)-1\bigr).
 \end{aligned}
\]
\end{claimproof}

We now prove part~(ii).
Put
\[
 k=\lceil r/h\rceil,\qquad q=\lfloor r/k\rfloor,
 \qquad
 B=R(n;r,r-1,1),\qquad M=R(n;r,r-1,r-h).
\]
Since $1\le r-h\le r-2$, monotonicity gives $n\le B\le M$.
We first show that the denominator in (ii) is positive.
Claim~\ref{cl:colorgrouping}, with $h=r-1$, gives
$B\le R(n;2)$.
By the ordinary Ramsey product inequality and $R(n;1)=n$,
\[
 R(n;r)\ge R(n;3)
 \ge(n-1)\bigl(R(n;2)-1\bigr)+1
 >R(n;2)\ge B,
\]
where the strict inequality follows from $n\ge3$ and
$R(n;2)>1$.
Thus $f(r-1)-f(1)>0$.

Suppose first that $q\ge2$.
Since $M\ge n$ and $(n-1)^q>n$,
\[
 (M-1)^q
 \ge M^q\left(\frac{n-1}{n}\right)^q
 >\frac{M^q}{n^{q-1}}.
\]
Claim~\ref{prop:terminalproduct} and $B\ge n$ therefore yield
\[
 R(n;r)\ge(M-1)^q+1
 >\frac{M^q}{n^{q-1}}
 \ge\frac{M^q}{B^{q-1}}.
\]
Taking logarithms gives
$\log(R(n;r)/B)>q\log(M/B)$, and hence
\[
 \frac{f(r-1)-f(r-h)}{f(r-1)-f(1)}
 =1-\frac{\log(M/B)}{\log(R(n;r)/B)}
 >1-\frac1q.
\]
If $q=1$, the corresponding non-strict bound follows from
$f(r-1)\ge f(r-h)$.

Since $h<r$, we have $k\ge2$.
The definition of $k$ gives $r\ge(k-1)h+1$, so
\[
 \frac rk\ge h-\frac{h-1}{k}\ge\frac{h+1}{2},
 \qquad
 q\ge\left\lfloor\frac{h+1}{2}\right\rfloor
   =\left\lceil\frac h2\right\rceil.
\]
Using the telescoping sums, we conclude that
\[
 \frac{\sum_{t=r-h}^{r-2}\Delta_t}
      {\sum_{t=1}^{r-2}\Delta_t}
 =\frac{f(r-1)-f(r-h)}{f(r-1)-f(1)}
 \ge1-\frac1q
 \ge1-\frac1{\lceil h/2\rceil}.
\]
This proves (ii).

For (iii), fix $r\ge4$ and put
\[
 m=\lfloor r/2\rfloor,\qquad
 a=R(n;r,r-1,r-2),\qquad b=R(n;m).
\]
Then $m\ge2$ and $a,b\ge n$.
The second inequality in Claim~\ref{prop:terminalproduct} gives
\[
 \frac{R(n;r)}a
 \ge\left(1-\frac1a\right)(b-1).
\]
Taking logarithms and dividing by $n$, we obtain
\[
 \begin{aligned}
 \Delta_{r-2}
 &\ge\frac{\log b}{n}
   +\frac{\log(1-1/a)+\log(1-1/b)}{n}\\
 &=\frac{\log b}{n}+o(1),
 \end{aligned}
\]
where the error tends to zero because $a,b\ge n$.
If $m=2$, Claim~\ref{prop:firstmoment} gives
\[
 \liminf_{n\to\infty}\frac1n\log R(n;2)
 \ge\frac12\log2.
\]
If $m\ge3$, Brada\v{c}~\cite[Theorem~1.6]{Bradac2026} gives
\[
 \log R(n;m)\ge\frac{m-1}{2}n\log2-O_m(1).
\]
Consequently,
\[
 \liminf_{n\to\infty}\Delta_{r-2}
 \ge\frac{m-1}{2}\log2
 =\frac{\lfloor r/2\rfloor-1}{2}\log2>0,
\]
proving (iii).

For (iv), let $r,n\ge3$ and put $a=R(n;r,r-1,r-2)$.
We show that $R(n;r)>a^2/n$ whenever $r$ is even or $n=3$.
If $r$ is even, Claim~\ref{prop:terminalproduct} with $h=2$
gives
\[
 R(n;r)\ge(a-1)^2+1\ge\frac{a^2}{2}>\frac{a^2}{n},
\]
since $(a-1)^2+1-a^2/2=(a-2)^2/2\ge0$ and $n\ge3$.

If $n=3$, Proposition~\ref{prop:exactinitialspectrum}(iii) gives
\[
 a=2\lceil r/2\rceil+1\le r+2.
\]
Repeated application of the ordinary Ramsey product inequality,
using $R(3;1)=3$, gives $R(3;r)\ge2^r+1$.
For $r\ge3$,
\[
 3(2^r+1)>(r+2)^2.
\]
Indeed, the difference equals $2$ at $r=3$ and increases
by $3\cdot2^r-2r-5>0$ when $r$ is replaced by $r+1$.
Thus
\[
 R(3;r)\ge2^r+1>\frac{(r+2)^2}{3}\ge\frac{a^2}{3}.
\]

In either case, $R(n;r)>a^2/n$.
Since $B=R(n;r,r-1,1)\ge n$, it follows that
\[
 \Delta_{r-2}
 =\frac1n\log\frac{R(n;r)}a
 >\frac1n\log\frac an
 \ge\frac1n\log\frac aB
 =f(r-2)-f(1)
 =\sum_{t=1}^{r-3}\Delta_t.
\]
This proves (iv) and completes the proof.
\end{proof}

\section{Elementary finite values}\label{sec:finite}

For completeness, we restate some elementary finite results
in the notation of this paper.
The triangle case and the cases $j=1,2$ of part~(ii)
correspond to classical observations discussed by Erd\H{o}s
and Gy\'arf\'as~\cite{EG1997}.
We give a self-contained proof, treating
$1\le j\le\lfloor n/4\rfloor$ uniformly.

\begin{proposition}\label{prop:exactinitialspectrum}
Let $n\ge3$ and $r\ge2$ be integers.
Then the following statements hold.
\begin{enumerate}
\item[(i)]
For every $1\le t<r$,
\(
 R(n;r,r-1,t)=n
 \;\Longleftrightarrow\;
 r-t\ge\binom n2.
\)
Consequently, if $r\ge\binom n2+1$, then
\(
 \Delta_t=0
 \;\left(1\le t<r-\binom n2\right),
 \;
 \sum_{t=r-\binom n2}^{r-2}\Delta_t=f(r-1)-f(1).
\)
Thus all the increase occurs in the final $\binom n2-1$
increments.

\item[(ii)]
For every integer $j$ with $1\le j\le\lfloor n/4\rfloor$,
if $r\ge\binom n2-j+1$, then
\(
 R\left(n;r,r-1,r-\binom n2+j\right)
 =
 1+\left\lfloor
      \frac{1+\sqrt{1+8(r+j-1)}}{2}
    \right\rfloor.
\)

In particular, for $n\ge4$ and $r\ge\binom n2$,
\(
 R\left(n;r,r-1,r-\binom n2+1\right)
 =
 1+\left\lfloor
      \frac{1+\sqrt{1+8r}}{2}
    \right\rfloor.
\)

\item[(iii)]
For $r\ge3$,
\(
 R(3;r,r-1,t)=
 \begin{cases}
  3, & 1\le t\le r-3,\\[2pt]
  2\lceil r/2\rceil+1, & t=r-2,\\[2pt]
  R(3;r), & t=r-1.
 \end{cases}
\)
\end{enumerate}
\end{proposition}

\begin{proof}
Put $m=\binom n2$.
We use the omitted-color formulation throughout:
$R(n;r,r-1,t)$ is the least $N$ such that every ordinary
coloring of $K_N$ with colors from $[r]$ contains an
$n$-vertex clique using at most $h=r-t$ colors.
We call a coloring avoiding if every $n$-vertex clique
uses at least $h+1$ colors.
In particular, $R(n;r,r-1,t)\ge n$.

For (i), suppose first that $h\ge m$.
Every coloring of $K_n$ uses at most $m$ colors, so
$R(n;r,r-1,t)=n$.
Conversely, if $h<m$, then $h+1\le\min\{m,r\}$.
We can therefore color the edges of $K_n$ using exactly
$h+1$ colors.
Its only $n$-vertex clique uses more than $h$ colors,
and hence $R(n;r,r-1,t)>n$.

If $r\ge m+1$, it follows that
$f(t)=\log n/n$ for $1\le t\le r-m$.
Thus $\Delta_t=0$ for $1\le t<r-m$, and telescoping gives
\[
 \sum_{t=r-m}^{r-2}\Delta_t
 =f(r-1)-f(r-m)
 =f(r-1)-f(1).
\]
The summation contains exactly $m-1$ increments,
proving (i).

For (ii), fix $j$ and $r$ as in the statement, and put
\[
 h=m-j,\qquad t=r-m+j.
\]
The assumptions imply $1\le h\le r-1$ and $1\le t<r$.
Fix an ordinary coloring $\psi:E(K_N)\to[r]$ with $N\ge n$.
For $E'\subseteq E(K_N)$, define
\[
 D(E')
 =
 |E'|-\bigl|\{\psi(e):e\in E'\}\bigr|.
\]
Adding an edge increases the number of edges by one
and the number of distinct colors by at most one.
Thus $D$ is nondecreasing under inclusion.

We claim that the coloring is avoiding if and only if
\[
 D(E(K_N))\le j-1.
\]
Suppose first that this inequality holds.
For every $n$-element vertex set $U$, monotonicity gives
\[
 D(E(K_N[U]))\le j-1.
\]
Since $K_N[U]$ has $m$ edges, it uses at least
$m-j+1=h+1$ colors.
Thus the coloring is avoiding.

Conversely, suppose that $D(E(K_N))\ge j$.
Choose one representative edge from each nonempty color class.
The number of edges outside this collection is exactly
$D(E(K_N))$, so we may select $j$ such edges.
Let $k\le j$ be the number of colors appearing on these
selected edges, and add the representative of each of
these $k$ colors.
The resulting edge set $E'$ has $j+k$ edges and exactly
$k$ colors, because the added representatives are distinct
from the selected edges and introduce no new colors.
Therefore
\[
 |E'|=j+k\le2j,
 \qquad
 D(E')=j.
\]
The edges in $E'$ have at most $4j\le n$ endpoints.
Since $N\ge n$, these endpoints are contained in an
$n$-element vertex set $U$.
By monotonicity,
\[
 D(E(K_N[U]))\ge D(E')=j.
\]
Hence $K_N[U]$ uses at most $m-j=h$ colors,
so the coloring is not avoiding.
This proves the claim.

Since at most $r$ colors are available, an avoiding coloring
must satisfy
\[
 \binom N2-r\le D(E(K_N))\le j-1,
\]
and consequently $\binom N2\le r+j-1$.
Conversely, if $\binom N2\le r+j-1$, color the edges using
exactly $\min\{r,\binom N2\}$ colors.
Then
\[
 D(E(K_N))
 =\max\left\{\binom N2-r,0\right\}
 \le j-1,
\]
so the claim shows that this coloring is avoiding.

Thus, for $N\ge n$, an avoiding coloring exists if and only if
$\binom N2\le r+j-1$.
The largest integer satisfying this inequality is
\[
 N_0=
 \left\lfloor
   \frac{1+\sqrt{1+8(r+j-1)}}{2}
 \right\rfloor.
\]
The assumption $r\ge m-j+1$ ensures that $N_0\ge n$.
An avoiding coloring therefore exists on $N_0$ vertices,
whereas every coloring on $N_0+1$ vertices contains
the required clique.
Consequently,
\[
 R(n;r,r-1,r-m+j)=N_0+1.
\]
Taking $j=1$ gives the particular case.

Finally, consider (iii).
For $n=3$, part~(i) gives
$R(3;r,r-1,t)=3$ whenever $1\le t\le r-3$.
At $t=r-1$, the required triangle is monochromatic, so
$R(3;r,r-1,r-1)=R(3;r)$.

It remains to consider $t=r-2$.
An avoiding coloring must assign three distinct colors
to every triangle.
This is equivalent to a proper edge-coloring of $K_N$,
in which any two edges sharing a vertex have different colors.
Indeed, two equally colored incident edges lie in a triangle
using at most two colors; conversely, the three edges of
a triangle are pairwise incident.

For $N\ge3$, the minimum number of colors in a proper
edge-coloring of $K_N$ is
\[
 \chi'(K_N)=
 \begin{cases}
  N-1, & N\text{ even},\\
  N,   & N\text{ odd}.
 \end{cases}
\]
We include the standard construction.
Suppose first that $N=2q$ is even, and label the vertices by
$\{\infty\}\cup\mathbb Z_{2q-1}$.
For each $a\in\mathbb Z_{2q-1}$, define
\[
 M_a=
 \bigl\{\{\infty,a\}\bigr\}
 \,\cup\,
 \bigl\{\{a+i,a-i\}:1\le i\le q-1\bigr\},
\]
with arithmetic modulo $2q-1$.
The residues $a+i$ and $a-i$, for $1\le i\le q-1$,
exhaust all residues other than $a$.
Thus each $M_a$ is a perfect matching.

Every edge $\{\infty,a\}$ belongs only to $M_a$.
For an edge $\{x,y\}$ with finite endpoints, the congruence
\[
 2a\equiv x+y\pmod{2q-1}
\]
has a unique solution, since $2q-1$ is odd.
This is precisely the index of the matching containing
$\{x,y\}$.
Hence the $N-1$ matchings $M_a$ partition all edges,
giving a proper edge-coloring with $N-1$ colors.
At least $N-1$ colors are necessary because the $N-1$
edges incident with any vertex must have distinct colors.

If $N$ is odd, every color class is a matching containing
at most $(N-1)/2$ edges.
Thus at least
\[
 \frac{\binom N2}{(N-1)/2}=N
\]
colors are necessary.
Deleting one vertex from a proper edge-coloring of $K_{N+1}$
gives a proper edge-coloring of $K_N$ with at most $N$ colors.
This proves the formula for $\chi'(K_N)$.

It follows that the largest order of a complete graph
admitting a proper edge-coloring with at most $r$ colors
is $r$ when $r$ is even and $r+1$ when $r$ is odd,
namely $2\lceil r/2\rceil$.
Therefore
\[
 R(3;r,r-1,r-2)=2\lceil r/2\rceil+1,
\]
which completes the proof.
\end{proof}

\section*{Declaration of Generative AI and AI-assisted technologies in the writing process}
During the preparation of this manuscript, the authors used generative AI to assist with English revision, organization and presentation of proofs, \LaTeX\ editing, and checking of mathematical arguments and references.
The authors reviewed and revised all AI-assisted content, independently verified all mathematical statements and cited references, and take full responsibility for the final manuscript.

\end{document}